\documentclass[10pt,a4paper]{article}
\makeatletter
\usepackage[dvipsnames]{xcolor}
\usepackage{lipsum}
\usepackage{tipa}
\usepackage{tikz}
\usepackage{tikz-cd}
\usepackage{float}
\usepackage{wrapfig}
\usetikzlibrary{matrix,arrows,decorations.pathmorphing}
\newcounter{magicrownumbers}

\usepackage[left=3cm,top=3cm,bottom=3cm,right=3cm,head=1cm,foot=1cm]{geometry}
\usepackage[pdftex,bookmarks=true,bookmarksnumbered=true]{hyperref}
\hypersetup{linkbordercolor=ForestGreen!25, citebordercolor=Goldenrod!60}
\DeclareMathAlphabet{\mathpzc}{OT1}{pzc}{m}{it}
\usepackage{graphicx}
\usepackage{stmaryrd}
\usepackage{twoopt}
\usepackage{amssymb}
\usepackage{amsmath}
\usepackage{amsthm}
\usepackage{mathabx} 
\usepackage{amsfonts}
\usepackage{amscd}
\usepackage{url}
\usepackage[all]{xy}

\newcommand{\lb}{\vspace{0,15cm} \noindent}

\newcommand{\Tr}{\mathrm{Tr}}

\newcommand{\hil}{\mathcal{H}}

\newcommand{\g}{\mathfrak{g}}

\newcommand{\cl}{\mathrm{cl}}

\newcommand{\5}{\hspace{0,5cm}}
\newcommand{\3}{\vspace{0,3cm}}

\title{Quantum Orbifolds}
\date{}
\author{Antti J. Harju\footnote{University of Regensburg} \footnote{harjuaj@gmail.com}}
\begin{document}
\maketitle

%
%
%
	
\begin{abstract}
This is a study of orbifold-quotients of quantum groups (quantum orbifolds $\Theta \rightrightarrows G_q$). These structures have been studied extensively in the case of the quantum $SU_2$ group. I will introduce a generalized mechanism which allows one to construct quantum orbifolds from any compact simple and simply connected quantum group. Associated with a quantum orbifold there is an invariant subalgebra as well as a crossed product algebra. For each spin quantum orbifold, there is a unitary equivalence class of Dirac spectral triples over the invariant subalgebra, and for each effective spin quantum orbifold associated with a finite group action, there is a unitary equivalence class of Dirac spectral triples over the crossed product algebra. As an application I will study a Hopf-equivariant Fredholm index problem.  \3



\end{abstract}

\section*{Introduction}

\noindent \textbf{I.1. Quantum Orbifolds.} There are well known cases of actions of the circle group $\mathbb{T}$ on the quantum group $SU_{2q}$, \cite{BF12}. The algebra of polynomial functions $\mathbb{C}[SU_{2q}]$ accepts left and right actions of the dual Hopf algebra $U_q(\mathfrak{su}_2)$ (Drinfeld-Jimbo algebra) which is a noncocommutative Hopf algebra. It is therefore rather surprising that the $\mathbb{T}$-algebra structures in $\mathbb{C}[SU_{2q}]$ exist. However, the key point is that the circle group $\mathbb{T}$ can be embedded onto a subset of group-like elements in the algebra of symmetries of $\mathbb{C}[SU_{2q}]$ and this makes the $\mathbb{T}$-algebra structures in $\mathbb{C}[SU_{2q}]$ possible. Another interesting feature in these actions is that, at least in certain special cases, the $\mathbb{T}$-algebra would seem to describe a noncommutative orbifold rather than a manifold. The orbifold phenomenon occurs naturally when the left and the right action are applied simultaneously.

In this work I will generalize the ideas discussed in the previous paragraph to the compact simple and simply connected quantum groups. The motivation is to construct models for potentially singular noncommutative spaces, rather than develop free actions. A maximal torus of the Lie group $G$ can be embedded onto a subset of group-like elements of the symmetry algebra of the quantum group. Then the torus, or any of its closed subgroup, is given an action on the quantum group  which applies the left and the right symmetries simultaneously. It is not difficult to construct such actions explicitly. In Section 3 I will study the case of quantum $SU_3$ in details and write down all the possibilities. Associated with any of these actions there are two algebras of interest: the crossed product algebra and the invariant subalgebra. However, the purpose of this work is to keep the discussion on a theoretical level and I will not make an attempt to classify all these algebras. I also want to point out, which should be evident by 2.2, that even in the case of quantum $SU_2$ there are undiscovered algebras arising from this theory. 

The quantum group $G_q$ together with the group action will be referred to a 'quantum orbifold'. However, one should be careful with this terminology: since we are applying the torus group $T$, the resulting isotropy groups on the classical geometric level are not always finite. They are finite if we choose a finite subgroup, and they are finite in some cases of the full torus action.

\3 \noindent \textbf{I.2. On Spectral Geometry of Noncommutative Orbifolds.} A fundamental problem is to develop a geometric theory for the quantum orbifolds. In differential geometry, an orbifold can be modeled in terms of the classical orbifold theory ("V-manifolds" \cite{Sat56}, \cite{Sat57}), or as orbifold groupoids in Lie groupoid theory\cite{MM03, MPro97}. In any case, the fundamental idea is to have a collection of locally defined smooth base manifolds subject to group actions and rules for gluing these together so that the orbit space of the system has the desired structure. There are many ways of how to realize the orbit space geometrically: for example, the teardrop orbifold which is one of the best known examples, can be modeled as a global quotient of $SU_2$ under an action of the circle group, or as a system of two 2-dimensional open discs glued together so that one is subject to an action of a finite rotation resulting in a cone shaped orbit space. One needs to set an equivalence relation between any models which describe the same orbit space, including the isotropy types of its singularities. This is what Morita equivalence does. For a deeper analysis of quantum orbifolds, one should have an understanding of the spectral triples arising from quantum orbifolds as well as a precise axiomatic formulation of the Morita equivalence in the theory of noncommutative orbifolds. In this work  the structure of the spectral triples will be studied.

Given a classical orbifold groupoid there are two relevant complex algebras which can be applied for a construction of a spectral triple. The convolution algebra of the orbifold groupoid is useful in the study of homology and K-theory. In particular one can apply the spectral triple theory for the convolution algebras in the study of Fredholm index problems, \cite{GL13}. Another approach is to build a spectral triple on the subalgebra of invariant smooth functions. The invariant subalgebra would seem to be the correct algebra for the study of metric properties since the geodesic lengths between points in the same orbit is considered to be zero.  However, in this case a difficulty arises since the universal differential complex built from the algebra of smooth invariant functions is not sufficiently large, \cite{RV08}. More precisely, if the orbifold is singular, one does not get the whole complex of invariant differential forms on the orbifold. 

In the quantum orbifold theory the crossed product algebra replaces the convolution algebra. I will construct Dirac spectral triples over the crossed product algebras and the invariant subalgebras under some natural conditions which are introduced in Section 4. The spectral triples developed in Section 4 are true geometric spectral triples which are built on the same principles as on the geometric level: Theorem 2 is a noncommutative analogue of Theorem 1 in \cite{Har14a} and Theorem 4 is a noncommutative analogue of Theorem 2 in \cite{Har14a}. In all of these 4 cases, the point is to construct an equivariant spin lift and a Dirac operator acting on the spinors which is invariant under the action of the orbifold. In the case of noncommutative orbifolds the construction is really only a matter of careful design: since the orbifold acts under its embedding onto the natural quantum group symmetries, one can use the $U_q(\g)$-invariant Dirac operator on a quantum group, \cite{NT10}, \cite{NT11b}. The spin condition introduced in 4.2 is the most general condition one can set to make the spectral triple well defined. In the case of $SU_{2q}$ spectral triples based on similar consideration have been developed in \cite{Har14b} and \cite{SV13}.

The primary open problem is to understand the Morita equivalence on the level of spectral triples of noncommutative orbifolds. This problem has been discussed in the case of geometric orbifold groupoids in \cite{Har14c}. In this case a Morita equivalence operates as a unitary equivalence on the invariant spectral triple over a compact spin orbifold. If the orbifold is also effective, then the unitary equivalence extends to define a transformation of spectral triples over the convolution algebra. The latter transformation is not a unitary equivalence in general. Based on this geometric study, one should be able to find algebraic axioms for the Morita transformations of spectral triples. The geometric study also suggest that an orbifold should be considered as nonsingular if and only if the convolution spectral triple is Morita equivalent to the invariant spectral triple. This should be viewed as a noncommutative geometric analogue of the geometric fact that an orbifold groupoid is Morita equivalent to the unit groupoid over its orbit space if and only if the orbifold is nonsingular. The freeness of noncommutative orbifolds has been recently discussed from the viewpoint of the invariant algebra in \cite{Bre14}, \cite{BS14}.  

As an application I will study the following Fredholm index problem. An invariant spectral triple can be coupled to an equivariant projective module which gives rise to a new invariant spectral triple. The index problem for the coupled Dirac operator will be studied. Notice that the origin of this problem is somewhere in the middle ground between the invariant formalism and the crossed product (equivariant) formalism. The computation of the index can be done with the standard mechanism in noncommutative geometry which in this case leads to a pairing of an equivariant K-theory and an invariant cyclic cohomology. However, here we have got to the situation where the universal differential complex of the invariant spectral triple is not large enough but the problem is fixed by introducing a larger complex of noncommutative invariant forms. In a general case, the universal complex is not sufficient because the equivariant K-theory carries some information which is not contained in the K-theory of the invariant subalgebra. Notice that the K-theory of the invariant algebra would map to the universal complex of the invariant algebra through the usual Chern character map, \cite{Con94}. So, one may view the obstruction to send the equivariant K-theory to the universal complex of the invariant algebra as a consequence of singularities on a noncommutative level: if the theory was free of singularities, then the equivariant K-theory should be isomorphic to the K-theory of the invariant algebra. 

\3 I acknowledge funding from the projects "Gerbes, K-theory, and Infinite-Dimensional Lie Groups" at the University of Helsinki and "SFB 1085 Higher Invariants" at the University of Regensburg. 

\3 \noindent \textbf{I.3. Notation.} Let $G$ be a compact simple and simply connected Lie group, and $\g$ the complexified Lie algebra. Fix a maximal torus $T$ in $G$ and denote by $\mathfrak{h}$ the corresponding Cartan subalgebra of $\g$. Choose a system of simple roots $\{\alpha_1, \ldots,  \alpha_n\}$. Let $a$ denote the $n \times n$-Cartan matrix and let $d_1, \ldots, d_n$ denote a set of coprime integers so that $(d_i a_{ij})_{ij}$ is symmetric. The bilinear form on $\mathfrak{h}^*$ is fixed so that $(\alpha_i, \alpha_j) = d_i a_{ij}$.  The generators $h_i$, $1 \leq i \leq n$ for the Cartan subalgebra are fixed so that $\alpha_j(h_i) = a_{ij}$. Let $P\subset \mathfrak{h}^*$ denote the set of integral weights and $P_+$ the subset of dominant integral weights. 

The symbol $q$ denotes a real number in $(0,1)$.

\section{Quantum Group Preliminaries}

\noindent \textbf{1.1.} Fix a compact simple and simply connected Lie group $G$ of rank $n$ and follow the notation of I.3. The Drinfeld-Jimbo algebra, $U_q(\mathfrak{g})$, is the polynomial algebra over $\mathbb{C}$ generated by $e_i, f_i, k_i, k_i^{- 1}$ for $ 1 \leq i \leq n$ which are subject to the relations 
\begin{eqnarray*}
& &[k_{i}, k_{j}] = 0,\5  k_i k_i^{-1} = 1 \5 k_{i} e_j k_{i}^{-1} = q_i^{\frac{a_{ij}}{2}} e_j,\5 k_{i} f_j k_{i}^{-1} = q_i^{-\frac{a_{ij}}{2}} f_j, \\
& &[e_i, f_j] = \delta_{ij}\frac{k^2_{i} - k_{i}^{-2}}{q_i -q_i^{-1}}, \\
& & \sum_{k=0}^{1-a_{ij}} (-1)^k {1 - a_{ij} \brack k}_{q_i} e_i^k e_j e_i^{1-a_{ij}-k} = 0 \5 \sum_{k=0}^{1-a_{ij}} (-1)^k {1 - a_{ij} \brack k}_{q_i} f_i^k f_j f_i^{1-a_{ij}-k} = 0
\end{eqnarray*}
where $q_i = q^{d_i}$ and
\begin{eqnarray*}
&& {m \brack k}_{q_i} = \frac{[m]_{q_i}!}{[k]_{q_i}![m-k]_{q_i}!},\hspace{0,5cm} [m]_{q_i} ! = [m]_{q_i} \cdots [2]_{q_i} [1]_{q_i}, \hspace{0,5cm} [n]_{q_i} = \frac{q_i^{n} - q_i^{-n}}{q_i - q_i^{-1}}.
\end{eqnarray*}
We fix the following Hopf structure in $U_q(\mathfrak{g})$:
\begin{eqnarray*}
\triangle_q(k_i) =k_i \otimes k_i\ & \triangle_q(e_i) = e_i \otimes k_i + k_i^{-1} \otimes e_i\  &\triangle_q(f_i) = f_i \otimes k_i + k_i^{-1} \otimes f_i  \\
S_q(k_i) = k_i^{-1} & S_q(e_i) = -q e_i & S_q(f_i) = -q^{-1}f_i \\
 \epsilon_q(k_i) = 1 & \epsilon_q(e_i) =  \epsilon_q(f_i) = 0&
\end{eqnarray*}
There is a $*$-structure in $U_q(\mathfrak{g})$, an anti-automorphism given on the generators by $e_i^* = f_i$, $f_i^* = e_i$ and $k_i^* = k_i$, and extended anti-linearly. Now $U_q(\g)$ is a Hopf $*$-algebra.  

Suppose that $(\varrho_q, M_q)$ is a $U_q(\g)$-representation. A vector $v \in M_q$ is called a weight $\mu$-vector if there is $\mu \in P$ such that $\varrho_q(k_i) v = q^{\mu(h_i)/2} v$. Denote by $M_q(\mu)$ the subspace of $M_q$ spanned by the weight $\mu$-vectors. A representation $(\varrho_{q}, M_q)$ is admissible if it has a weight decomposition: $M_q = \bigoplus_{\mu \in P} M_q(\mu)$. A representation $(\varrho_{q}, M_q)$ is called a highest weight representation if there are vectors $v_{\lambda} \in V$ and $\lambda \in \mathfrak{h}^*$ such that
\begin{eqnarray*}
\varrho_{\lambda,q}(k_i) v_{\lambda} = q_i^{\frac{\lambda(h_i)}{2}} v_{\lambda}, \5 \varrho_{\lambda,q}(e_i) v_{\lambda} = 0
\end{eqnarray*} 
for each $i$. The isomorphism classes of irreducible highest weight representations are indexed by the dominant integral weights, $\lambda \in P_+$. Consider the category of finite dimensional admissible $U_q(\g)$-representations $(\mathfrak{C},\g,q)$. This category is semisimple. The simple objects are the irreducible highest weight representations. There is a well known braided monoidal structure in the category $(\mathfrak{C},\g,q)$, see e.g. \cite{NT11}. 

\3 \noindent \textbf{1.2.} We fix a highest weight $*$-representation $\varrho_{\lambda q}: U_q(\g) \rightarrow \mathcal{B}(M_{\lambda})$ for each  $\lambda \in P_+$. The matrix elements of the $*$-representations satisfy $\varrho_{\lambda q}(x^*) = (\varrho_{\lambda q}(x))^{\dagger}$ for all $x \in U_q(\g)$. Define the dual compact quantum group of $G_q$ by 
\begin{eqnarray*}
\widehat{\mathbb{C}[G_q]} = \bigoplus_{\lambda \in P_+} \mathcal{B}(M_{\lambda,q}).
\end{eqnarray*}
There is a $*$-bialgebra structure in $\widehat{\mathbb{C}[G_q]}$ which is inherited from $U_q(\g)$ since the representations $\varrho_{\lambda q}: U_q(\g) \rightarrow \mathcal{B}(M_{\lambda,q})$ are all surjective. The von Neumann algebra $W^*(G_q)$ is the $l_{\infty}$-product of the finite dimensional algebras $\mathcal{B}(M_{\lambda})$ for all $\lambda\in P_+$. This algebra should be understood as a completion of the dual quantum group.  We also define the algebraic product
\begin{eqnarray*}
\mathcal{U}(G_q) = \prod_{\lambda \in P_+} \mathcal{B}(M_{\lambda,q}).
\end{eqnarray*}
The algebra $U_q(\g)$ can also be realized as a subalgebra in $\mathcal{U}(G_q)$. There is a $*$-structure and a $*$-homomorphism $\triangle_q : W^*(G_q) \rightarrow W^*(G_q) \otimes W^*(G_q)$ which extends to a $*$-homomorphism:
\begin{eqnarray*}
\triangle_q: \mathcal{U}(G_q) \rightarrow \mathcal{U}(G_q \times G_q) = \prod_{(\lambda, \mu) \in P_+ \times P_+} \mathcal{B}(M_{\lambda,q} \otimes M_{\mu,q}). 
\end{eqnarray*}
The universal $R$-matrix $\mathcal{R} \in \mathcal{U}(G_q \times G_q)$ is a unique element satisfying $\triangle_q^{\text{op}}(\cdot) = \mathcal{R} \triangle_q(\cdot) \mathcal{R}^{-1}$ and if $(\varrho_{\lambda q}, M_{\lambda})$ is a finite dimensional representation with a highest weight vector $v_{\lambda, \lambda}$, and if $(\varrho_{\lambda' q}, M_{\lambda'})$ is a finite dimensional representation with a lowest weigth vector $v_{\lambda', \mu}$, then 
\begin{eqnarray*}
(\varrho_{\lambda q} \otimes \varrho_{\lambda' q})(\mathcal{R})(v_{\lambda, \lambda} \otimes v_{\lambda', \mu}) = q^{(\lambda, \mu)} v_{\lambda, \lambda} \otimes v_{\lambda', \mu}. 
\end{eqnarray*}

The simple objects on the category $(\mathfrak{C}, \g,q)$ and in the category of finite dimensional admissible $\g$-representations are parametrized by the same set, $P_+$.  The dimensions of the simple objects match. So, the algebras $W^*(G_q)$ and $W^*(G)$ are isomorphic. The group von Neumann algebra $W^*(G)$ is defined in terms of the representations of $\g$ and the bialgebraic structures in $U(\g)$ determine bialgebraic structures in $W^*(G)$. The essential differences between $W^*(G)$ and $W^*(G_q)$ are described in the following theorem which we recall from \cite{NT10}.

\3 \noindent \textbf{Theorem 1.} There exists a $*$-homomorphism $\phi: W^*(G_q) \rightarrow W^*(G)$ which identifies the centers, and a unitary $\mathcal{F} \in W^*(G) \otimes W^*(G)$ so that 
\begin{quote}
\textbf{1.} $(\phi \otimes \phi) \triangle_q = \mathcal{F} \triangle(\cdot) \mathcal{F}^{-1}$, 

\textbf{2.} $(\varepsilon \otimes \iota)(\mathcal{F}) = (\iota \otimes \varepsilon)(\mathcal{F}) = 1$, 

\textbf{3.} $(\phi \otimes \phi)(\mathcal{R}) = \mathcal{F}_{21}q^t \mathcal{F}^{-1}$,

\textbf{4.} $\Phi = (\iota \otimes \triangle)(\mathcal{F}^{-1})(1 \otimes \mathcal{F}^{-1})(\mathcal{F} \otimes 1)(\triangle \otimes \iota)(\mathcal{F})$ defines an associativity morphism in the braided tensor category of admissible $\g$-modules with the braiding determined by the symmetric group. This associator coincides with the Drinfeld's associator.  
\end{quote}

\noindent In the previous theorem, we used $t = \sum_k x_k \otimes x_k \in \g \otimes \g$ where the sum is over an index set of a basis of $\g$ which is normalized so that $(x_k ,x_l) = - \delta_{kl}$ with respect to the Killing form. The element $\mathcal{F}$ of Theorem 1 will be referred to a twist. The isomorphism of Theorem 1 extends to the isomorphism 
\begin{eqnarray*}
\phi: \mathcal{U}(G_q) \rightarrow \mathcal{U}(G). 
\end{eqnarray*}
which identifies the centers.

\3 \noindent \textbf{1.3.} Define the vector space of matrix elements: 
\begin{eqnarray*}
\mathbb{C}[G_q] = \bigoplus_{\lambda \in P_+} M^*_{\lambda,q} \otimes M_{\lambda,q}.
\end{eqnarray*}
Consider the standard nondegenerate bilinear pairing $\mathbb{C}[G_q] \otimes U_q(\g) \rightarrow \mathbb{C}$ determined by 
\begin{eqnarray*}
t^{\lambda}_{\mu \nu}(x) = u^*_{\lambda, \mu}(\varrho_{\lambda q}(x) u_{\lambda, \nu}).
\end{eqnarray*}
There is a product in $\mathbb{C}[G_q]$ which is determined by the relations 
\begin{eqnarray*}
t_1 t_2(x) = t_1(x') t_2(x'')
\end{eqnarray*}
for all $t_1, t_2 \in \mathbb{C}[G_q]$ and $x \in U_q(\g)$. This product makes $\mathbb{C}[G_q]$ a unital algebra: the unit is given by a matrix coefficient of the trivial representation. The vector space basis in $\mathbb{C}[G_q]$ can be chosen by 
\begin{eqnarray*}
t^{\lambda}_{\mu \nu} = u^*_{\lambda, \mu} \otimes u_{\lambda, \nu}. 
\end{eqnarray*}
where $\lambda$ runs over $P_+$, the vectors $u_{\lambda, \nu}$ define a weight-basis of the module $M_{\lambda}$ and $\mu$, $\nu$ run over all the parameters of this basis. The Hopf structure in $\mathbb{C}[G_q]$ can be fixed by requiring that $U_q(\g)$ and $\mathbb{C}[G_q]$ are Hopf dual with respect to the pairing. So, the coproduct is given by $\triangle_q(t^{\lambda}_{\mu \nu}) = \sum_{\eta} t^{\lambda}_{\mu \eta} \otimes t^{\lambda}_{\eta \nu}$. The dual $*$-structure in $\mathbb{C}[G_q]$ is defined by $t^*(x) = \overline{t(S(x)^*)}$.  

The enveloping algebra $U_q(\mathfrak{\g})$ acts on $\mathbb{C}[SU_{2,q}]$ under the right regular representation by
\begin{eqnarray*}
\partial_q(x) t^{\lambda}_{\mu \nu} = u^*_{\lambda, \mu} \otimes \varrho_{\lambda q}(x) u_{\lambda, \nu}
\end{eqnarray*}
Consider the following algebra automorphism in $U_q(\g)$: 
\begin{eqnarray*}
\vartheta(k_i) = k_i^{-1}, \5 \vartheta(e_i) = - f_i, \5 \vartheta(f_i) = - e_i. 
\end{eqnarray*}
It is straightforward to check that 
\begin{eqnarray*}
\triangle_q(S_q(\vartheta(x))) = S_q(\vartheta(x')) \otimes S_q(\vartheta(x'')). 
\end{eqnarray*}
Then we define the left action of $U_q(\g)$ on $\mathbb{C}[G_q]$ by 
\begin{eqnarray*}
l_q(x) t^{\lambda}_{\mu \nu} = \varrho^*_{\lambda,q}(\vartheta(x))u^*_{\lambda, \mu} \otimes u_{\lambda, \nu}.
\end{eqnarray*}

Associated to the Haar state in $\mathbb{C}[G_q]$, there is the GNS triple $(L^2(G_q), \pi, \xi)$ where $\xi$ is the cyclic vector of the representation. Let $\eta: \mathbb{C}[G_q] \rightarrow L^2(G_q)$ be the canonical map $\eta: t \mapsto t \xi$. The basis is given by $\eta(t^{\lambda}_{\mu \nu})$ where $\lambda$ runs over $P_+$ and $\mu,\nu$ over all their possible values. The representation of $\mathbb{C}[G_q]$ on $L^2(G_q)$ is defined by 
\begin{eqnarray*}
\pi(t^{\lambda}_{\mu \nu}) \eta(t^{\lambda'}_{\mu' \nu'})   = \eta(t^{\lambda}_{\mu \nu} t^{\lambda'}_{\mu' \nu'}),
\end{eqnarray*}
and extended linearly. The algebra $U_q(\g)$ has the commuting right and left actions on $L^2(G_q)$ which are given on the basis by 
\begin{eqnarray*}
\partial_q(x) \eta(t^{\lambda}_{\mu \nu}) = \eta(\partial_q(x)t^{\lambda}_{\mu \nu}),\5 l_q(x) \eta(t^{\lambda}_{\mu \nu}) = \eta(l_q(x)t^{\lambda}_{\mu \nu}) 
\end{eqnarray*}
for all $x \in U_q(\g)$ and $t^{\lambda}_{\mu \nu} \in \mathbb{C}[G_q]$.

\3 \noindent \textbf{Proposition 1.} The GNS representation is left and right equivariant with respect to the action of $U_q(\g)$ in the sense that 
\begin{eqnarray*}
\partial_q(x)(\pi_h(t^{\lambda}_{\mu \nu}) \eta(t^{\lambda'}_{\mu' \nu'})) &=& \pi_h(\partial_q(x') t^{\lambda}_{\mu \nu}) \eta(\partial_q(x'')t^{\lambda'}_{\mu' \nu'}) \\
l_q(x)(\pi_h(t^{\lambda}_{\mu \nu}) \eta(t^{\lambda'}_{\mu' \nu'})) &=& \pi_h(l_q(x') t^{\lambda}_{\mu \nu}) \eta(l_q(x'')t^{\lambda'}_{\mu' \nu'})
\end{eqnarray*}
for all $x \in U_q(\g)$.

\3 \noindent Proof. For the first identity, notice that the left side can be written by $\eta (\partial_q(x)(t^{\lambda}_{\mu \nu} t^{\lambda'}_{\mu' \nu'}))$ whereas the right 
side can be written by $\eta(\partial_q(x') t^{\lambda}_{\mu \nu} \partial_q(x'')t^{\lambda'}_{\mu' \nu'})$. Since $\eta$ is injective, it is sufficient to show that   
\begin{eqnarray}\label{prop1}
 \partial_q(x)(t^{\lambda}_{\mu \nu} t^{\lambda'}_{\mu' \nu'}) = \partial_q(x') t^{\lambda}_{\mu \nu} \partial_q(x'')t^{\lambda'}_{\mu' \nu'}
\end{eqnarray}
By the definition of the product, we have $t^{\lambda}_{\mu \nu} t^{\lambda'}_{\mu' \nu'}(x) = t^{\lambda}_{\mu \nu}(x') t^{\lambda'}_{\mu' \nu'}(x'')$ for all indices. This means that the matrix elements of the right representation $\partial_q$ are the same on both sides of \eqref{prop1}, and so the equality holds. 

For the equivariance under the left representation, note that 
\begin{eqnarray*}
t^{\lambda}_{\mu \nu} t^{\lambda'}_{\mu' \nu'}(S_q \vartheta(x)) &=& t^{\lambda}_{\mu \nu}(S_q(\vartheta (x))') t^{\lambda'}_{\mu' \nu'}(S_q(\vartheta (x))'') \\
&=& t^{\lambda}_{\mu \nu}(S_q(\vartheta (x'))) t^{\lambda'}_{\mu' \nu'}(S_q(\vartheta (x''))).
\end{eqnarray*}
So, the claim follows from the same consideration as in the right case. \5 $\square$ \3

The left and right $U_q(\g)$ actions above determine two commuting actions of $W^*(G_q)$ on $L^2(G_q)$ which are still denoted by $\partial_q$ and $l_q$. These actions extend to give $\mathcal{U}(G_q)$ two commuting actions by densely defined unbounded operators on $L^2(G_q)$. The equivariance properties of Proposition 1 continue to hold for these actions.

\section{The Algebras of Quantum Orbifolds}
  
\noindent \textbf{2.1.} Let $T$ be a maximal torus of $G$. The maximal torus may be viewed as a subset of the set of group-like elements in $\mathcal{U}(G_q)$ with respect to the homomorphism $\triangle_q$. Concretely, the basis of the Cartan subalgebra $\mathfrak{h}$ in $\mathcal{U}(G_q)$ is obtained by applying $\text{log}_q$ on the vectors $k_1, \ldots, k_n$, and then the torus is constructed by exponentiation. Denote by $\mathcal{U}(T)$ the image of the torus in $\mathcal{U}(G_q)$. Suppose that $K$ is a closed subgroup of $T$. Consider any action $\Theta: K \times \mathbb{C}[G_q] \rightarrow \mathbb{C}[G_q]$ of $K$ on the space $\mathbb{C}[G_q]$ which factors as a product of a left and a right actions, 
\begin{eqnarray}\label{def}
\Theta(\sigma, t) = l_q(\sigma_{(1)}) \partial_q(\sigma_{(2)})t
\end{eqnarray}
for all $t \in \mathbb{C}[G_q]$ and $\sigma \in K$, so that $\sigma \mapsto [\sigma_{(1)}, \sigma_{(2)}]$  is a map $K \rightarrow (\mathcal{U}(T) \times \mathcal{U}(T)) / \sim$ and the equivalence relations is defined by 
\begin{eqnarray*}
(\sigma_1,\sigma_2) \sim (\sigma_1 \mu, \mu^{-1} \sigma_2)
\end{eqnarray*}
if $\sigma_1, \sigma_2 \in \mathcal{U}(T)$ and $\mu \in \mathcal{U}(T) \cap Z \mathcal{U}(G_q)$ where $Z\mathcal{U}(G_q)$ is the center of $\mathcal{U}(G_q)$. The most fundamental case is that $\sigma \mapsto \sigma_{(1)}$, $\sigma \mapsto \sigma_{(2)}$ is a pair of group homomorphisms $K \rightarrow \mathcal{U}(T)$ and they are given an action on $\mathbb{C}[G_q]$ as in \eqref{def}. However, this is not the most general case, as will be evident from the examples in 2.3 and Section 3 where all these actions are found in the case of $SU_{3q}$. 

Let us write $\sigma \triangleright t$ to denote $\Theta(\sigma,t)$ in $\mathbb{C}[G_q]$. A quantum group $G_q$ equipped with an action $\Theta$ satisfying the conditions  introduced above is referred to a quantum orbifold $\Theta \rightrightarrows G_q$. 

\3 \noindent \textbf{Proposition 2.} If $\Theta \rightrightarrows G_q$ is a quantum orbifold associated with an action of the group $K$, then $\mathbb{C}[G_q]$ is a $K$-algebra. 

\3 \noindent Proof. The antihomomorphism $S_q \circ \vartheta$ in the left action $l_q$ restricts in $\mathcal{U}(T)$ as the identity homomorphism. Observe that the action \eqref{def} is independent on the choice of the equivalence class of $[\sigma_{(1)}, \sigma_{(2)}]$ for all $\sigma \in K$. We need to show that the action $\Theta$ determines a map $K \rightarrow \text{Aut}(\mathbb{C}[G_q])$. It is sufficient to show  that the action is compatible with the product: if $t_1$ and $t_2$ are in $\mathbb{C}[G_q]$ and $\sigma$ is in $K$ then $\sigma \triangleright t_1 t_2 = (\sigma \triangleright t_1) (\sigma \triangleright t_2)$. The following holds for all $x \in \mathcal{U}(G_q)$ 
\begin{eqnarray*}
(\sigma \triangleright t_1 t_2)(x) &=& t_1 t_2(S_q(\vartheta(\sigma_{(1)})) x \sigma_{(2)}) \\
&=& t_1 t_2(\sigma_{(1)} x \sigma_{(2)}) \\
&=& t_1((\sigma_{(1)} x \sigma_{(2)})') t_2((\sigma_{(1)} x \sigma_{(2)})'') \\
&=& t_1(\sigma_{(1)} x' \sigma_{(2)}) t_2(\sigma_{(1)} x'' \sigma_{(2)}) \\
&=& (\sigma \triangleright t_1)(x') (\sigma \triangleright t_2)(x'') 
\end{eqnarray*}
where we have used the property that the images of $\sigma_{(i)}$ are grouplike. \5 $\square$\3

The invariant quantum orbifold algebra $\mathbb{C}[G_q]^{\Theta}$ consists of the vectors that are invariant under the action $\Theta$. 

\3 \noindent \textbf{Proposition 3.} The invariant subalgebra $\mathbb{C}[G_q]^{\Theta}$ is an involutive algebra. 

\3 \noindent Proof. The dual matrix element of $t^{\lambda}_{\mu \nu}$ satisfies 
\begin{eqnarray*}
(t^{\lambda}_{\mu \nu})^*(x) &=& \overline{t^{\lambda}_{\mu \nu}(S_q(x)^*)} \\
&=& \overline{v^*_{\lambda \mu} (\varrho_{\lambda q}(S_q(x)^*) v_{\lambda \nu}) } \\
&=&  \overline{v^*_{\lambda \mu} (\varrho_{\lambda q}(S_q(x))^{\dagger} v_{\lambda \nu}) } \\
&=&  v^*_{\lambda \mu} (\varrho_{\lambda q}(S_q(x))^{t} v_{\lambda \nu}) \\
&=& t^{\lambda^{\vee}}_{-\mu, -\nu}(x)
\end{eqnarray*}
where $t^{\lambda^{\vee}}_{-\mu, -\nu}$ are  matrix elements of the dual representation $(\varrho^*_{\lambda,q}, M_{\lambda^{\vee}})$ corresponding to the weights $- \mu$ and $- \nu$. So, given an action $\Theta$ as in \eqref{def}, and $\sigma \in K$, we can find $y_{(1)}$ and $y_{(2)}$ in $\mathfrak{h}$ so that $\sigma_{(a)} = e^{i y_{(a)}}$ for $a = 1,2$ and then  
\begin{eqnarray}\label{star1}
\sigma \triangleright t^{\lambda}_{\mu \nu} = e^{\frac{i}{2}(\mu(y_{(1)}) + \nu(y_{(2)}))}t^{\lambda}_{\mu \nu}
\end{eqnarray}
and, according to the formula above, we have 
\begin{eqnarray}\label{star2}
\sigma \triangleright (t^{\lambda}_{\mu \nu})^* &=& \sigma \triangleright t^{\lambda^{\vee}}_{-\mu,- \nu} \\
&=& e^{-\frac{i}{2}(\mu(y_{(1)}) + \nu(y_{(2)}))} (t^{\lambda}_{\mu \nu})^* \nonumber
\end{eqnarray}
So, the invariant subalgebra is stable under the $*$-operation. \5 $\square$  \3

It is often convenient to view a quantum group $\mathbb{C}[G_q]$ as a $*$-algebra defined in terms of generators and relations. In this case, it is sufficient to determine the action $\Theta$ on the generators and  apply Proposition 2 to define the action on the whole algebra. Then one can use \eqref{star1} and \eqref{star2} to define the action on the $*$-conjugates.  

Suppose that $K$ is a finite subgroup of $T$. The crossed product quantum orbifold algebra $\mathbb{C}[K \ltimes G_q]$ is the linear space of all finite formal combinations of elements in $K$ with coefficients in $\mathbb{C}[G_q]$, and the product is defined by 
\begin{eqnarray*}
(\sigma_1\cdot t_1) (\sigma_2\cdot t_2) = (\sigma_1 \sigma_2 \cdot t_1(\sigma_1 \triangleright t_2))
\end{eqnarray*}
for all $\sigma_1, \sigma_2 \in K$ and $t_1, t_2 \in \mathbb{C}[G_q]$, extended linearly. This is a $*$-algebra with 
\begin{eqnarray*}
(\sigma \cdot t)^* = \sigma^{-1}\cdot(\sigma^{-1} \triangleright t^*) 
\end{eqnarray*}
for all $\sigma \in K$ and $t \in \mathbb{C}[G_q]$, extended anti-linearly. 

One could also define the crossed product algebra in the general case where $K$ is not assumed to be finite but this requires some modifications. The algebra $\mathbb{C}[G_q]$ needs to be completed to a topological algebra and then the elements of the crossed product algebra are continuous maps from $K$ to the algebra. In particular one can consider the $C^*$-algebraic completion of  $\mathbb{C}[G_q]$ which leads to a $C^*$-crossed product algebra.

\3 \noindent \textbf{2.2.} Let us first consider $G = SU_2$ very briefly. Take the maximal torus to be the diagonal subgroup and identify it with $\mathbb{T}$. The generator $h$ of the Cartan subalgebra of $\g$ is chosen to be $h = \text{diag}(1,-1)$. The highest weights $\lambda$ are parametrized by the half-integers $\frac{1}{2} \mathbb{N}_0$ (the highest weight vector $v_{\lambda}$ in $M_{\lambda}$ corresponds to $\lambda \in \frac{1}{2} \mathbb{N}_0$ if $\varrho_{\lambda q}(k) v_{\lambda} = q^{\lambda} v_{\lambda}$).  Following the notation of \cite{Har14b} we denote by $t^{\frac{1}{2}}_{mn}$ the matrix elements of the defining representations, $m,n \in \{ \pm \frac{1}{2} \}$. We use the symbols $\alpha$ and $\beta$ for the generators of the $*$-algebra $\mathbb{C}[SU_{2q}]$ which can be chosen by
\begin{eqnarray*}
\alpha = t^{\frac{1}{2}}_{\frac{1}{2},\frac{1}{2}} \5 \text{and} \5 \beta = t^{\frac{1}{2}}_{\frac{1}{2},-\frac{1}{2}}. 
\end{eqnarray*}
The algebraic relations $\beta \alpha = q\alpha \beta$, $\beta^* \alpha = q \alpha \beta^*$, $\beta \beta^* = \beta^* \beta$, $\alpha \alpha^* + \beta \beta^* = 1$ and $\alpha^* \alpha + q^2 \beta^* \beta = 1$ are straightforward to derive using the standard $*$-representatin of $U_q(\g)$ and the Hopf-duality with $U_q(\mathfrak{su}_2)$, see \cite{Har14b}.  

Suppose that $k$ and $l$ are positive coprime integers and $K = \mathbb{T}$. Consider the action $\Theta$ which factors through the map $K \rightarrow (\mathcal{U}(T) \times \mathcal{U}(T)) / \sim $ given by 
\begin{eqnarray}\label{22}
e^{i \phi} \mapsto [e^{i\phi(l-k) h/2}, e^{-i\phi(l+k) h/2}]. 
\end{eqnarray}
The matrix elements transform according to 
\begin{eqnarray*}
e^{i \phi} \triangleright \alpha = e^{- i k \phi} \alpha, \5  e^{i  \phi} \triangleright \beta = e^{+i l \phi} \beta
\end{eqnarray*}
which defines a $\mathbb{T}$-algebra structure in $\mathbb{C}[SU_{2q}]$. The fixed point algebra is known as the quantum weighted projective plane $\mathbb{C}[\mathbb{WP}_{klq}]$. If $k = l = 1$ we have the standard quantum sphere. If $k = 1$ and $l \in \mathbb{N}$ we have a quantum teardrop. Now take $k = 1$, $l \in \mathbb{N}$ and $K = \mathbb{Z}_p$ for some prime number $p$. Then $\phi$ gets values in $\{\frac{2 \pi j}{p} : j = 1, \ldots, p \}$ and the associated fixed point algebra is the quantum lens space algebra, \cite{BF12}.

\section{Quantum Orbifolds on $SU_{3q}$}
 
\textbf{3.1.} The algebra $U_q(\mathfrak{su}_3)$ is generated by $e_i, f_i, k_i, k_i^{-1}$ with $i = 1,2$. The highest weights of $SU_3$ can be parametrized by the set $\mathbb{N}_0  \times \mathbb{N}_0$. The trivial representation corresponds to $(0,0)$ and $(1,0)$ corresponds to the defining representation. The dual of $M_{(1,0)}$ is isomorphic to $M_{(0,1)}$. We shall also need to apply the 6-dimensional representation associated with the highest weight $(2,0)$. Let us fix the following notation
\begin{eqnarray*}
\lambda_0 = (0,0),\5 \lambda_1 = (1,0),\5 \lambda_1^{\vee} = (0,1), \5 \lambda_2 = (2,0).
\end{eqnarray*}
We fix $(\varrho_{\lambda_1 q}, M_{\lambda_1})$ to be the 3-dimensional $*$-representation
\begin{eqnarray*}
\varrho_{\lambda_1 q}(e_1) &=& e_{12}, \\
 \varrho_{\lambda_1 q}(e_2) &=& e_{23},\\ 
\varrho_{\lambda_1 q}(k_1) &=& \text{diag}(q^{\frac{1}{2}}, q^{-\frac{1}{2}}, 1),\\
 \varrho_{\lambda_1 q}(k_2)& =& \text{diag}(1, q^{\frac{1}{2}}, q^{-\frac{1}{2}})
\end{eqnarray*}
where the basis vectors are parametrized by $\{1,2,3\}$ and we use $e_{ij}$ the matrix with only nonzero entry $1$ on the $i$'th row an $j$'th column. The $*$-representation $(\varrho_{\lambda_1^{\vee} q}, M_{\lambda_1^{\vee}})$ is fixed by
\begin{eqnarray*}
\varrho_{\lambda_1^{\vee} q}(e_1) &=& e_{21},\\ 
\varrho_{\lambda_1^{\vee} q}(e_2) &=& e_{32},\\
\varrho_{\lambda_1^{\vee} q}(k_1) &=& \text{diag}(q^{- \frac{1}{2}}, q^{\frac{1}{2}}, 1),\\ 
\varrho_{\lambda_1^{\vee} q}(k_2) &=& \text{diag}(1, q^{- \frac{1}{2}}, q^{\frac{1}{2}}).
\end{eqnarray*}
In this module, the parameter $3$ corresponds to the highest weight vector.  For $M_{\lambda_2}$ we use the $*$-representation 
\begin{eqnarray*}
&\varrho_{\lambda_2 q}(e_1)& = \sqrt{[2]} e_{12} + \sqrt{[2]} e_{23} + e_{45}\\
&\varrho_{\lambda_2 q}(e_2)& = e_{24} + \sqrt{[2]} e_{35} + \sqrt{[2]} e_{56}\\
&\varrho_{\lambda_2 q}(k_1)& = q e_{11} + e_{22} + q^{-1} e_{33} + q^{\frac{1}{2}} e_{44} + q^{-\frac{1}{2}} e_{55} + e_{66}.\\
&\varrho_{\lambda_2 q}(k_2)& = e_{11} + q^{\frac{1}{2}}e_{22} + q^{1} e_{33} + q^{-\frac{1}{2}} e_{44} + e_{55} + q^{-1} e_{66}. 
\end{eqnarray*}
The basis vectors are parametrized by $\{1,2,3,4,5,6\}$ and the highest weight vector corresponds to the parameter 1. 

The generators of the algebra $\mathbb{C}[SU_{3q}]$ can be taken to be any 9 linearly independent matrix elements of $M_{\lambda_1}^* \otimes M_{\lambda_1}$. We fix the generators $t_{ij} = v_i^* \otimes v_j$ for all $1 \leq i,j \leq 3$, where $v_i$ are the basis vectors of $M_{\lambda_1}$ for which the representation of $U_q(\mathfrak{su}_3)$ is given by $\varrho_{\lambda_1 q}$. The commutation relations for these generators can be solved by studying the tensor product decomposition
\begin{eqnarray*}
M_{\lambda_1} \otimes M_{\lambda_1} \simeq M_{\lambda_2} \oplus M_{\lambda_1^{\vee}}. 
\end{eqnarray*}
We can use orthogonal Clebsch-Gordan matrices to change the tensor product basis on the left side, on which $U_q(\mathfrak{su}_3)$ acts rhtough the coproduct, to the basis of the module on the right hand side. The product rule reads
\begin{eqnarray}\label{product2}
t_{mn} t_{m' n'} = \sum_{\kappa \in \{\lambda_1^{\vee}, \lambda_2 \}} \begin{bmatrix} \lambda_1 & \lambda_1 & \kappa \\ m & m' & k_{\kappa} \end{bmatrix}  \begin{bmatrix} \lambda_1 & \lambda_1 & \kappa \\ n & n' & k_{\kappa}' \end{bmatrix} t^{\kappa}_{k_{\kappa},k_{\kappa}'}. 
\end{eqnarray}  
This can be derived exactly as in the case of $SU_2$, see e.g. \cite{Har14b}, and can be specialized for any products in $\mathbb{C}[SU_{3q}]$. Relative to the representations fixed in the previous paragraph, the nonzero Clebsch-Gordan coefficients are 
\begin{eqnarray*}
&& \begin{bmatrix}
\lambda_1 & \lambda_1 & \lambda_1^* \\ 1 & 2 & 1
\end{bmatrix} = 
\begin{bmatrix}
\lambda_1 & \lambda_1 & \lambda_1^* \\ 1 & 3 & 2
\end{bmatrix} =
\begin{bmatrix}
\lambda_1 & \lambda_1 & \lambda_1^* \\ 2 & 3 & 3
\end{bmatrix} = \sqrt{\frac{q}{[2]}}  \\
&&  \begin{bmatrix}
\lambda_1 & \lambda_1 & \lambda_1^* \\ 2 & 1 & 1
\end{bmatrix} = 
\begin{bmatrix}
\lambda_1 & \lambda_1 & \lambda_1^* \\ 3 & 1 & 2
\end{bmatrix} =
\begin{bmatrix}
\lambda_1 & \lambda_1 & \lambda_1^* \\ 3 & 2 & 3
\end{bmatrix} = - \sqrt{\frac{1}{q[2]}}.
\end{eqnarray*}
and
\begin{eqnarray*}
&&\begin{bmatrix} \lambda_1 & \lambda_1 & \lambda_2 \\ 1 & 1 & 1 \end{bmatrix} = 
\begin{bmatrix} \lambda_1 & \lambda_1 & \lambda_2 \\ 2 & 2 & 3 \end{bmatrix} = 
\begin{bmatrix} \lambda_1 & \lambda_1 & \lambda_2 \\ 3 & 3 & 6 \end{bmatrix} = 1 \\
&&\begin{bmatrix} \lambda_1 & \lambda_1 & \lambda_2 \\ 2 & 1 & 2 \end{bmatrix} = 
\begin{bmatrix} \lambda_1 & \lambda_1 & \lambda_2 \\ 3 & 1 & 4 \end{bmatrix} = 
\begin{bmatrix} \lambda_1 & \lambda_1 & \lambda_2 \\ 3 & 2 & 5 \end{bmatrix} = \sqrt{\frac{q}{[2]}} \\
&&\begin{bmatrix} \lambda_1 & \lambda_1 & \lambda_2 \\ 1 & 2 & 2 \end{bmatrix} = 
\begin{bmatrix} \lambda_1 & \lambda_1 & \lambda_2 \\ 1 & 3 & 4 \end{bmatrix} = 
\begin{bmatrix} \lambda_1 & \lambda_1 & \lambda_2 \\ 2 & 3 & 5 \end{bmatrix} = \sqrt{\frac{1}{q[2]}}.
\end{eqnarray*}
The highest weight vectors of the tensor product modules are normalized so that these matrices are orthogonal. After a straightforward but slightly tedious computations one finds the following formulas
\begin{eqnarray*}
&&t_{ij} t_{ik} = q t_{ik} t_{ij} \5 \text{if} \5  j < k \\
&&t_{ik} t_{jk} = q t_{jk} t_{ik} \5 \text{if} \5  i < j \\
&&t_{ij} t_{kl} = t_{kl} t_{ij} \5 \text{if} \5 i < k, j > l \\
&&[t_{ij}, t_{kl}]  =(q -q^{-1}) t_{il} t_{kj} \5  \text{if} \5 i < k \5 \text{and} \5 j < l.
\end{eqnarray*}
The unit is a matrix element of the trivial representation. The module $M_{\lambda_1} \otimes M_{\lambda_1^{\vee}}$ is isomorphic to a sum of the adjoint representation and the trivial representation. A direct computations with the Clebsch-Gordan matrices gives 
\begin{eqnarray*}
1 = t^*_{33} t_{33} + q^{-2} t_{32}^* t_{32} + q^{-4} t^*_{31} t_{31}. 
\end{eqnarray*}
The $*$-conjugates are given in the next proposition.

It is straightforward to write down formulas for the coproduct and the counit. The antipode sends the matrix elements in the representation $M_{\lambda_1}$ to the matrix elements in its dual $M_{\lambda_1^{\vee}}$. Explicit formulas can be obtained by studying the Hopf duality which gives the following. 

\3 \noindent \textbf{Proposition 4.} The antipode and the $*$-structure in the algebra $\mathbb{C}[G_q]$ are fixed by the following relations: 
\begin{eqnarray*}
S_q(t_{11}) = t^{\vee}_{11} \5 & S_q(t_{12}) = -q t^{\vee}_{21} & \5 S_q(t_{21}) = -q^{-1} t^{\vee}_{12} \\
S_q(t_{22}) = t^{\vee}_{22} \5 & S_q(t_{23}) = -q t^{\vee}_{32} & \5 S_q(t_{32}) = -q^{-1} t^{\vee}_{23} \\
S_q(t_{33}) = t^{\vee}_{33} \5 & S_q(t_{13}) = q^{2} t^{\vee}_{31} & \5 S_q(t_{31}) = q^{-2} t^{\vee}_{13}.
\end{eqnarray*}
and $t^*_{ij} = S_q(t_{ji})$ for all $1 \leq i,j \leq 3$. \3

In the operator algebraic literature, see e.g. \cite{Wor87}, an antipode of a generator is typically written as a sum of quadratics of the generators. However, note that $M_{\lambda_1^{\vee}}$ is a submodule of $M_{\lambda_1} \otimes M_{\lambda_1}$ and so we can always write the elements $t^{\vee}_{ij}$ as quadratics of the generators using the rule \eqref{product2} and the Clebsch-Gordan coefficients above. 

\3 \noindent \textbf{3.2.} Let $T$ be the maximal torus of $G$ and $K$ a closed subgroup. Let $h_1$ and $h_2$ be the generators of the Cartan subalgebra. There is a group isomorphism in $T$ which sends the pair $(T,K)$ to $(\mathbb{T}^2, K_1 \times K_2)$ so that $K_i$ is $\mathbb{T}$ or a finite cyclic subgroup of roots of unity for $i = 1,2$. We wish to find all the actions of $K$ on $\mathbb{C}[SU_{3q}]$ satisfying the conditions of 2.1. However, it is sufficient to solve this problem for the groups $K_1$ and $K_2$ independently since the action of $K$ is fully determined by the actions of $K_1$ and $K_2$ and the actions of $K_1$ and $K_2$ can be defined independently on each other.  So, let us take $K_2 = 0$ (the one element subgroup of $\mathbb{T}$) and let $K_1$ be any closed subgroup of $\mathbb{T}$. Observe that any action $\Theta$ factor through the map $K \rightarrow (\mathcal{U}(T) \times \mathcal{U}(T)) / \sim$ given by 
\begin{eqnarray}\label{32}
(e^{i \phi}, 1) \mapsto [e^{i \phi y_{(1)}^1 h_1} e^{i \phi y_{(1)}^2 h_2} , e^{i \phi y_{(2)}^1 h_1} e^{i \phi y_{(2)}^2 h_2}]
\end{eqnarray} 
for some real numbers $y_{(a)}^b$, $a,b \in \{1,2\}$. Notice that the inner functions need to be  linear in $\phi$, otherwise the action could not be well defined. We can fix the generators of $\mathfrak{h}$ by $h_1 = \text{diag}(1,-1,0)$ and $h_2 = \text{diag}(0,1,-1)$.
 
\3 \noindent \textbf{Proposition 5.} If $K_1$ is a closed subgroup of $\mathbb{T}$ and $K_2 = 0$, then $\Theta$ determined by \eqref{32} defines a quantum orbifold if and only if the real numbers $y_{(a)}^b$ are chosen by   
\begin{eqnarray*}
y_{(1)}^1 = k_1^1+  \frac{x}{3}  \5 y_{(1)}^2 = k_1^2+\frac{2x}{3}  \5 y_{(2)}^1 = k_2^1 +\frac{2x}{3}  \5 y_{(2)}^2 = k_2^2 + \frac{x}{3} 
\end{eqnarray*}
for some $x = 0,1,2$ and $k_{a}^b \in \mathbb{Z}$ with $a,b \in \{1,2\}$. 

\3 \noindent Proof. The only nontrivial part to check is the $2 \pi$-periodicity. Observe that 
\begin{eqnarray*}
\mathcal{U}(T) \cap Z \mathcal{U}(SU_{3q}) = \{1, e^{i\frac{2\pi}{3} h_1} e^{i\frac{4\pi}{3} h_2} ,  e^{i\frac{4\pi}{3} h_1} e^{i\frac{2\pi}{3} h_2} \}
\end{eqnarray*}
which can be identified with the group $\mathbb{Z}_3$. If the value of the map \eqref{32} is evaluated modulo the action of $\mathcal{U}(T) \cap Z \mathcal{U}(SU_{3q})$, then \eqref{32} is indeed $2 \pi$-periodic if and only if $x = 0,1,2$ and $k_{a}^b \in \mathbb{Z}$. \5 $\square$\3

It is a simple task to write down the action on the generators of $\mathbb{C}[SU_{3q}]$ with any choice of $x$ and $y_{(a)}^b$. Notice that if $K_1 = \mathbb{Z}_p$, then it is sufficient to let $k_{b}^a$ get values in $\{0, \ldots, p-1\}$ only. Moreover, the same analysis as above can be applied with  the subgroup $K_2$ which lies in the second cartesian component of $\mathbb{T}^2$, and the same list of available parameters applies. 

Consider the case $K = T \simeq \mathbb{T}^2$ acting on the group $SU_3$ under the adjoint action $T \times SU_3 \rightarrow SU_3$; $(\sigma, g) \mapsto \sigma g \sigma^{-1}$. This is an obvious case of a singular action since the maximal torus in $SU_3$ lies in the singular locus. In fact, this is a case of a singular action groupoid which is not a classical orbifold since the isotropy group along the maximal torus is the torus group $K = T$. However, if we take $K$ to be a finite subgroup and let it act under the adjoint action we get a classical orbifold model. In the quantum group theory, this corresponds to the action factoring through
\begin{eqnarray*}
(e^{i \phi},e^{i \theta}) \mapsto [e^{i \phi h_1} e^{i \theta h_2} , e^{-i \phi h_1} e^{-i \theta h_2}].
\end{eqnarray*}
The generators $t_{11}$, $t_{22}$ and $t_{33}$ are invariant and the rest of the generators transform as 
\begin{eqnarray*}
e^{i \phi} \triangleright t_{12} = e^{2 i \phi} t_{12}&  \5  & e^{i \theta} \triangleright t_{12} = e^{-i \theta} t_{12} \\
e^{i \phi} \triangleright t_{13} = e^{i \phi} t_{13}&  \5  & e^{i \theta} \triangleright t_{13} = e^{i \theta} t_{13} \\
e^{i \phi} \triangleright t_{21} = e^{-2i \phi} t_{21}&  \5  & e^{i \theta} \triangleright t_{21} = e^{i \theta} t_{21} \\
e^{i \phi} \triangleright t_{23} = e^{-i \phi} t_{23}&  \5  & e^{i \theta} \triangleright t_{23} = e^{2i \theta} t_{23} \\
e^{i \phi} \triangleright t_{31} = e^{-i \phi} t_{31}&  \5  & e^{i \theta} \triangleright t_{31} = e^{-i \theta} t_{31} \\
e^{i \phi} \triangleright t_{32} = e^{i \phi} t_{32}&  \5  & e^{i \theta} \triangleright t_{32} = e^{-2i \theta} t_{32}.
\end{eqnarray*}
The groups $K_1$ and $K_2$ in $\mathbb{T}^2$ can also be restricted to finite subgroups of $\mathbb{T}$, and the parameters $\phi$ and $\theta$ are defined accordingly.

\section{Spectral Geometry on Quantum Orbifolds}

\noindent \noindent \textbf{4.1.} Let $\text{cl}(\g)$ denote the complexified Clifford algebra associated with the Killing form in $\g$ and let $\gamma: \g \rightarrow \text{cl}(\g)$ be the canonical embedding satisfying $\gamma(x)^2 = (x,x)1$ for all $x \in \g$. There is a Lie algebra homomorphism 
\begin{eqnarray*}
\widetilde{\text{ad}}: \g \rightarrow \text{cl}(\g); \5 \widetilde{\text{ad}}(x) = \frac{1}{4} \sum_k \gamma(x_k) \gamma([x,x_k])
\end{eqnarray*}
where $x_k$ is a basis of $\g$ for which the relations $(x_k, x_j) = - \delta_{kj}$ hold. Then, 
\begin{eqnarray*}
\gamma([x,y]) = [\widetilde{\text{ad}}(x), \gamma(y)]
\end{eqnarray*}
for all $x, y \in \g$. Fix a complex irreducible spinor module (irreducible $\cl(\g)$-representation) $(s, \Sigma)$. Now $\g$ is represented on $\Sigma$ through $s \circ \widetilde{\text{ad}}$. We shall use the following classical Dirac operator 
\begin{eqnarray*}
\mathcal{D} = \sum_k \Big( x_k \otimes \gamma(x_k) + \frac{1}{2} \gamma(x_k) \widetilde{\text{ad}}(x_k) \Big)
\end{eqnarray*}
which is an element in $U(\g) \otimes \text{cl}(\g)$. Let $\partial$ denote the right regular representation of $U(\g)$ on $L^2(G)$. Then $D = (\partial \otimes s) \mathcal{D}$ is an unbounded densely defined self-adjoint operator on $L^2(G) \otimes \Sigma$.

Suppose that $\phi$ and $\mathcal{F}$ are fixed (recall Theorem 1). Following \cite{NT10} we can define a quantum group Dirac operator as an element in $\mathcal{U}(G_q) \otimes \text{cl}(\g)$ by
\begin{eqnarray*}
\mathcal{D}_q = (\phi^{-1} \otimes \iota) \Big( (\iota \otimes \widetilde{\text{ad}})(\mathcal{F}) \mathcal{D} (\iota \otimes \widetilde{\text{ad}})(\mathcal{F}^*) \Big). 
\end{eqnarray*}
Now $D_q = (\partial_q \otimes s) \mathcal{D}_q$ is an unbounded densely defined self-adjoint operator on $\hil = L^2(G_q) \otimes \Sigma$. The group $\mathcal{U}(T)$ acts on $\hil$ under the following bounded representations 
\begin{eqnarray*}
\sigma \mapsto l_q(\sigma) \otimes \iota  \5 \text{and} \5 \sigma \mapsto \partial_q(\sigma) \otimes s \circ \widetilde{\text{ad}} \circ \phi(\sigma)
\end{eqnarray*}
for all $\sigma \in \mathcal{U}(T)$ and the Dirac operator $D_q$ commutes with both actions. 

Let us write $s_q(\sigma) = s \circ \widetilde{\text{ad}} \circ \phi(\sigma)$ for all $\sigma \in \mathcal{U}(T)$. 

\3 \noindent \textbf{4.2.} If the dimension of $G$ is an odd integer, then there are two isomorphism classes of irreducible spinor modules and if the dimension is even, then there is one isomorphism class. In the latter case, the spinor module has a decomposition $\Sigma = \Sigma_+ \oplus \Sigma_-$ so that the operators of the representation $s \circ \widetilde{\text{ad}}$ are even operators with respect to this decomposition, and the components $\Sigma_+$ and $\Sigma_-$ are $\g$-modules. Moreover, there is the self-adjoint chirality operator $\omega$ and $\Sigma_{\pm}$ are the $\pm 1$ eigenspaces of $\omega$.

Consider a quantum orbifold $\Theta \rightrightarrows G_q$ associated with an action of a closed subgroup $K$ in $T$. As a bounded action, $\Theta$ extends to define an action on the Hilbert space $L^2(G_q)$. We wish to lift this action to the Hilbert space of spinors $\hil$. An action is called a lift of $\Theta$ if it factors according to 
\begin{eqnarray*}
\sigma \triangleright (\eta(t) \otimes v) = \eta(\sigma \triangleright t) \otimes (\hat{\varrho}(\sigma) v) \5 t \otimes v \in \mathbb{C}[G_q] \otimes \Sigma
\end{eqnarray*}
so that in the odd (respectively, even) dimensional case, $\hat{\varrho}$ is a representation $\hat{\varrho}: K \rightarrow \mathcal{B}(\Sigma)$ (respectively, $K \rightarrow \mathcal{B}(\Sigma_+) \oplus \mathcal{B}(\Sigma_-)$), and the assignments
\begin{eqnarray*}
\sigma \mapsto \hat{\varrho}(\sigma) (\cdot) \hat{\varrho}(\sigma)^{-1} \5 \text{and} \5 \sigma \mapsto s_q(\sigma_{(2)}) (\cdot) s_q(\sigma_{(2)})^{-1}
\end{eqnarray*}
define the same linear transformations in the space $\mathcal{B}(\Sigma)$ (respectively, $\mathcal{B}(\Sigma_+) \oplus \mathcal{B}(\Sigma_-)$) for all $\sigma \in K$. If a lift exists, then $\Theta \rightrightarrows G_q$ is called a spin quantum orbifold. The lifts are referred to spin structures. Notice that the spin structures depend on the choice of the $*$-isomorphism $\phi$ because $s_q$ does. 

\3 \noindent \textbf{Example 1.} Consider the case $G = SU_2$ and $T = \mathbb{T}$. The irreducible spinor modules are 2-dimensional and the action $s_q$ makes a spinor module an irreducible $U_q(\mathfrak{su}_2)$-module of highest weight $\frac{1}{2}$.  The quantum orbifolds associated with the action \eqref{22} in 2.2 are all spin. The different spin structures are studied in the reference \cite{Har14b}, and they are parametrized by $\mathbb{Z}$. The spin structures arise since any 1-dimensional representation of $\mathbb{T}$ can be mapped to the center of $\mathcal{B}(\Sigma)$.

\3 This example reveals the rather obvious point that in the case of a spin quantum orbifold, the spin structures come in great numbers. This is not geometrically realistic, however, the condition  for the spin lift is exactly what is needed for the constructions of spectral triples. The difference is that in a purely geometric level, the lifting problem is affiliated with a lifting problem of an equivariant $SO_n$-bundle but we do not need to worry about this.

\3 \noindent \textbf{Example 2.} Let $G = SU_3$. We follow the notation of 3.2. Recall that  all actions of $K_1$ and $K_2$ on $SU_{3q}$ are mutually commutative and independent on each other and so the spin condition can be studied independently. So, let us take $K_1 = \mathbb{T}$ and $K_2 = 0$. The basis of the irreducible spinor module can be chosen so that $\Sigma = \Sigma_+ \oplus \Sigma_-$ decomposes into two irreducible components $(s_{q \pm},\Sigma_{\pm})$ under the representation $s_q$ which are both isomorphic to the adjoint representation $M_{(1,1)}$. Moreover, the basis can be chosen so that the Cartan elements $h_1$ and $h_2$ (considered as elements in $\mathcal{U}(G_q))$ act on these modules by
\begin{eqnarray*}
s_{q \pm}(h_1) = \text{diag}(1,2,-1,0,0,1,-2,-1) \5 \text{and} \5 s_{q \pm}(h_2) = \text{diag}(1,-1,2,0,0,-2,1,-1) 
\end{eqnarray*}  
In the notation of 3.2, the action of the element $\sigma_{(2)}$ on the spinor module reads
\begin{eqnarray*}
s_{q \pm}(\sigma_{(2)}) &=& \text{diag} \Big(e^{i \phi(k^1_2 + k^2_2 + x)}, e^{i \phi(2 k^1_2 - k^2_2 + x)}, e^{i \phi(-k^1_2 + 2 k^2_2)},1, \\
             &&1, e^{i \phi(k_2^1 - 2 k_2^2)}, e^{i \phi(-2 k_2^1 + k_2^2 - x)}, e^{i \phi(-k_2^1 - k_2^2 - x)}\Big)
\end{eqnarray*}
for all $x = 0,1,2$. So, the map $\sigma \mapsto s_{q}(\sigma_{(2)})$ is a well defined reducible representation of $K$, and $\hat{\varrho}(\sigma) = s_q(\sigma_{(2)})$ defines one spin structure. Now the spin structures are parametrized by $\mathbb{Z} \times \mathbb{Z}$ and they arise from the representations of $\mathbb{T}$ whose matrix elements lie in the center of $\mathcal{B}(\Sigma_{\pm})$. It follows that the quantum orbifolds $\Theta \rightrightarrows SU_{3q}$ are all spin: if $K_1 = K_2 = \mathbb{T}$, then the map, $\sigma \mapsto s_{q}(\sigma_{(2)})$ gives a spin structure, and this spin structure restricts to all closed subgroups.   

\3 It is not obvious that every quantum groupoid is spin. In the second example we found that for $x = 1,2$ the map $\sigma \mapsto s_q(\sigma_{(2)})$ gives a representation. However, this observation is specific to the highest weight $(1,1)$: for example for the highest weights $(1,0)$ or $(0,1)$ this would not lead to a well defined representation for the values $x = 1,2$. So, one might view this as an accidental coincidence. 

Let $\Theta \rightrightarrows G_q$ be a spin quantum orbifold. As a consequence of Proposition 1 we have that 
\begin{eqnarray*}
\sigma \triangleright \pi(t) \psi = \pi(\sigma \triangleright t) (\sigma \triangleright \psi)
\end{eqnarray*}
for all $\psi \in \hil$, $t \in \mathbb{C}[G_q]$ and $\sigma \in K$. It follows that $\pi \otimes \iota$ restricts to define a representation of $\mathbb{C}[G_q]^{\Theta}$ on the invariant subspace $\hil^{\Theta}$. However, let us not write the second tensor component explicitly. 

\3 \noindent \textbf{Theorem 2.} Let $\Theta \rightrightarrows G_q$ be a spin quantum orbifold. The collection $(\mathbb{C}[G_q]^{\Theta}, D_q, \hil^{\Theta})$ defines a finite summable spectral triple which is even if dim($G$) is even. Up to a unitary equivalence of spectral triples, this model is independent on the choices of $\mathcal{F}$ and $\phi$. 

\3 \noindent Proof. The algebra $\mathbb{C}[G_q]^{\Theta}$ is a $*$-algebra by 2.1 and the representation $\pi$ is faithful on $\hil^{\Theta}$ because it is a restriction of a faithful representation on $\hil$. The quantum Dirac operator commutes with the action of $K$ on $\hil$. To see this let us study the commuting right and lend actions independently: 
\begin{eqnarray*}
& & (\partial_q(\sigma_{(2)}) \otimes \hat{\varrho}(\sigma)) D_q \\
&=&  (\partial_q(\sigma_{(2)}) \otimes \hat{\varrho}(\sigma)) D_q (\partial_q(\sigma_{(2)})^{-1} \otimes \hat{\varrho}(\sigma)^{-1})(\partial_q(\sigma_{(2)}) \otimes \hat{\varrho}(\sigma)) \\
&=& (\partial_q(\sigma_{(2)}) \otimes s_q(\sigma_{(2)})) D_q (\partial_q(\sigma_{(2)})^{-1} \otimes s_q(\sigma_{(2)})^{-1})(\partial_q(\sigma_{(2)}) \otimes \hat{\varrho}(\sigma)) \\
&=& D_q (\partial_q(\sigma_{(2)}) \otimes \hat{\varrho}(\sigma)) 
\end{eqnarray*}
where the last equality follows from the $\mathcal{U}(T)$-invariance of $\mathcal{D}_q$, and  
\begin{eqnarray*}
(l_q(\sigma_{(1)}) \otimes \iota) D_q = D_q (l_q(\sigma_{(1)}) \otimes \iota)
\end{eqnarray*}
for all $\sigma \in K$. It follows that $D_q$ restricts to define an unbounded densely defined self-adjoint operator on $\hil^{\Theta}$. The commutators $[D_q, \pi(t)]$ extend to bounded operators because this holds in the larger Hilbert space $\hil$, \cite{NT10} Theorem 3.7. Moreover, the summability of $D_q$ is finite since this is the case for $D_q$ in the larger space $\hil$.

Consider the case of even dimensional $G$. Since the action of $K$ on $\Sigma = \Sigma_+ \oplus \Sigma$ decomposes as $\Sigma_{\pm} \rightarrow \Sigma_{\pm}$, this action commutes with the chiral grading $\omega$ on the quantum group model over $\mathbb{C}[G_q]$. Thus, $\omega$ restricts on $\hil^{\Theta}$ to define a chiral grading. 

 If $\mathcal{F}$, $\mathcal{F}'$ is a pair of twists, then by \cite{NT11b}, Theorem 6.1.(i) the associated Dirac operators are exactly the same. Then suppose that $\phi$, $\phi'$ is a pair of $*$-isomorphisms $\mathcal{U}(G_q) \rightarrow \mathcal{U}(G)$. There is a unitary element $u$ so that $\phi'(\cdot) = u \phi(\cdot) u^*$. The Dirac operators are related as follows, \cite{NT11b} Theorem 6.1.(ii)
\begin{eqnarray*}
\mathcal{D}_q' = (1 \otimes \widetilde{\text{ad}}(u)) \mathcal{D}_q (1 \otimes \widetilde{\text{ad}}(u^*)).
\end{eqnarray*}
It remains to prove that the unitary transformation $(1 \otimes s \circ \widetilde{\text{ad}}(u))$ on $\hil$ restricts to define an equivalence of invariant spectral triples. The chirality operator commutes with $(1 \otimes s \circ \widetilde{\text{ad}}(u))$ because the image of $\widetilde{\text{ad}}$ is even with respect to the Clifford grading.  Suppose that $\sigma \in K$ and $t \otimes v \in \mathbb{C}[G_q] \otimes \Sigma_+$, then 
\begin{eqnarray*}
(1 \otimes s \circ \widetilde{\text{ad}}(u)) \sigma \triangleright (t \otimes v) &=& (\sigma \triangleright t) \otimes (s \circ \widetilde{\text{ad}}(u) \hat{\varrho}(\sigma)v) \\
&=& (\sigma \triangleright t) \otimes \zeta(s \circ \widetilde{\text{ad}}(u) s_{q }(\sigma_{(2)}) v) \\
&=& (\sigma \triangleright t) \otimes \zeta(s \circ \widetilde{\text{ad}}(u \phi(\sigma_{(2)}))v) \\ 
&=& (\sigma \triangleright t) \otimes \zeta(s \circ \widetilde{\text{ad}}(\phi'(\sigma_{(2)}) u)v) \\ 
&=& (\sigma \triangleright t) \otimes (\hat{\varrho}'(\sigma)s \circ \widetilde{\text{ad}}(u) v) \\
&=& \sigma \triangleright  (1 \otimes s \circ \widetilde{\text{ad}}(u))(t \otimes v)
\end{eqnarray*}
for some $\zeta$ in the center of $\mathcal{B}(\Sigma_+)$. The same analysis shows that the action of $K$ commutes with this unitary operator on $\mathbb{C}[G_q] \otimes \Sigma_-$. So $(1 \otimes s \circ \widetilde{\text{ad}}(u))$ sends the invariant subspaces to each other. Since the representation of $\mathbb{C}[G_q]^{\Theta}$ commute with this operator we get a unitary equivalence of spectral triples. \5 $\square$

\3 \noindent \textbf{4.3.} In what follows we see that the universal differential algebra associated with the spectral triple of Theorem 2 might not be large enough for the study of index problems with coupling determined by an equivariant module. For this reason we introduce an additional piece of data, the complex $C_{\bullet}$ where $C_k$ is the subspace of $\mathbb{C}[G_q]^{\otimes k+1}$ consists of the invariant elements under 
\begin{eqnarray*}
\sigma \triangleright (t_0 \otimes \cdots \otimes t_k) = (\sigma \triangleright t_0) \otimes \cdots \otimes (\sigma \triangleright t_k).
\end{eqnarray*}
for each $k \in \mathbb{N}_0$. So, in particular $C_0 = \mathbb{C}[G_q]^{\Theta}$. With each $c_k \in C_k$ we associate a universal differential form by 
\begin{eqnarray*}
c_k = t_0 \otimes \cdots \otimes t_k \mapsto \pi(t_0) [D_q, \pi(t_1)] \cdots [D_q, \pi(t_k)].
\end{eqnarray*}
By the analysis of spectral triples on quantum groups, \cite{NT10}, these forms extend to bounded operators on $\hil$. Since $D_q$ is $K$-invariant, the invariant forms restrict to define bounded operators on $\hil^{\Theta}$. 

\3 \noindent \textbf{4.4.} Suppose that $G$ is an even dimensional Lie group of dimension $2n$ and $K$ is a closed subgroup of $T$. Let $\mathbb{C}[K]$ be the group Hopf algebra of $K$: the Hopf structure is the primitive one. The action of $K$ on $\mathbb{C}[G_q]$ determines an action of $\mathbb{C}[K]$ on $\mathbb{C}[G_q]$ which we denote by $\Theta$. 

Following \cite{NT04}, we can identify the $\mathbb{C}[K]$-equivariant projective modules on $\mathbb{C}[G_q]$ with the $\mathbb{C}[K]$-invariant idempotents in $\text{end}(V) \otimes \mathbb{C}[G_q]$ where $V$ ranges over the finite representations $(\rho,V)$ of the Hopf algebra $\mathbb{C}[K]$. The invariance is under the action
\begin{eqnarray*}
\sigma \triangleright (A \otimes t) = \rho(\sigma) A \rho(\sigma)^{-1} \otimes (\sigma \triangleright t)
\end{eqnarray*}
for $A \otimes t \in \text{end}(V) \otimes \mathbb{C}[G_q]$. Two such elements $p \in \text{end}(V) \otimes \mathbb{C}[G_q]$ and $p' \in \text{end}(V') \otimes \mathbb{C}[G_q]$ are defined to be equivalent if there exists $\mathbb{C}[K]$-invariant elements $\gamma \in \text{hom}(V,V') \otimes \mathbb{C}[G_q]$ and $\gamma' \in \text{hom}(V',V) \otimes \mathbb{C}[G_q]$ so that $\gamma \gamma' = p'$ and $\gamma' \gamma = p$. One can take direct sums of these modules and define the monoid structure in the set of equivalence classes of  $\mathbb{C}[K]$-equivariant projective modules. The Hopf-equivariant (even) K-theory group $K^{\Theta}(\mathbb{C}[G_q])$ is the Grothendieck group completion. 

Consider a Hopf equivariant module $p \in \text{end}(V) \otimes \mathbb{C}[G_q]$ for some finite representation $V$ of the Hopf algebra $\mathbb{C}[K]$. Then, for all $k \in \mathbb{N}_0$ we define the components of the character of $p$ by 
\begin{eqnarray*}
\text{ch}^{\Theta}_{2k}(p) = \text{tr}(p \otimes \cdots \otimes p),
\end{eqnarray*}
where $p$ appears $2k + 1$ times on the right side and tr is the usual matrix trace operation applied for the tensor products of matrix valued algebra elements, \cite{Con94}. 

\3 \noindent \textbf{Proposition 6.} The characters of the Hopf-equivariant projective modules are elements in the complex $C_{\bullet}$: 
$\text{ch}^{\Theta}_{2k}(p) \in C_{2k}$ for all $k \in \mathbb{N}_0$. 

\3 \noindent Proof. Suppose that $\sigma \in \mathbb{C}[K]$ and $p = \sum_i A^i \otimes t_i$. Then 
\begin{eqnarray*}
\sigma \triangleright \text{ch}^{\Theta}_{2k}(p) &=& \sigma \triangleright \sum_{i_0, \ldots, i_{2k}} \text{tr}(A^{i_0}  \cdots A^{i_{2k}}) t_{i_0} \otimes \cdots \otimes t_{i_{2k}} \\
&=&  \sum_{i_0, \ldots, i_{2k}} \text{tr}(A^{i_0}  \cdots A^{i_{2k}}) (\sigma \triangleright t_{i_0}) \otimes \cdots \otimes (\sigma \triangleright t_{i_{2k}}) \\
&=& \sum_{i_0, \ldots, i_{2k}} \text{tr} (\rho(\sigma) A^{i_0} \rho(\sigma)^{-1} \cdots \rho(\sigma) A^{i_{2k}} \rho(\sigma)^{-1})(\sigma \triangleright t_{i_0}) \otimes \cdots \otimes (\sigma \triangleright t_{i_{2k}}) \\
&=& \text{ch}^{\Theta}_{2k}(\sigma \triangleright p) =  \text{ch}^{\Theta}_{2k}(p). 
\end{eqnarray*}
The last equality follows from the invariance of $p$. \5 $\square$ \3

Nothing in the previous proof indicates that the characters would be an element in the universal differential complex determined by the invariant subalgebra: the individual tensor components may not lie in $\mathbb{C}[G_q]^{\Theta}$ or $\mathbb{C}[G_q]^{\Theta}/\mathbb{C}$. It is well known in the classical differential geometry that in the case of nontrivial isotropy, the universal complex determined by the invariant smooth functions is too small to be a target for nontrivial character maps from equivariant K-theory, see \cite{Har14a}, \cite{RV08}. On the other hand, if the group action is free, then this complex is sufficient. Since we are modeling a noncommutative orbifold, we have introduced the larger complex $C_{\bullet}$ which was proven to be large enough. 

\3 Denote by $C^k$ the space of linear functional on $C_k$. Since the action of $\mathbb{C}[K]$ on $C_k$ applies the cocommutative coproduct, the cyclicity and Hochschild coboundary operators restrict to define the linear maps
\begin{eqnarray*}
\lambda: C^k \rightarrow C^k \5 \text{and} \5 b: C^k \rightarrow C^{k+1}
\end{eqnarray*}
for all $k \in \mathbb{N}_0$. The cyclic subcomplex of $C^{\bullet}$ consist of the cyclic functionals: $c \in C^{\bullet}$ such that $\lambda(c) = c$. The Hopf-invariant cyclic cohomology $HC^{\bullet}_{\Theta}(\mathbb{C}[G_q])$ is the cohomology of the cyclic subcomplex $C^{\bullet}$ associated with the coboundary operator $b$. 

\3 \noindent \textbf{Proposition 7.} 
\begin{quote}
\textbf{1.} For all $k \in \mathbb{N}_0$, there are the pairings 
\begin{eqnarray*}
&&K^{\Theta}(\mathbb{C}[G_q]) \times HC^{2k}_{\Theta}(\mathbb{C}[G_q]) \rightarrow \mathbb{C} \\
&&([p], [c]) \mapsto c(\text{ch}^{\Theta}_{2k}(p)). 
\end{eqnarray*}
\textbf{2.} If $2n$ is the dimension of $G$, then
\begin{eqnarray*}
\text{ch}^{2n}_{\Theta}(t_0 \otimes \cdots \otimes t_{2n}) = \frac{(-1)^n}{2} \Tr\ ( \omega F_q [F_q, \pi(t_0)] \cdots  [F_q, \pi(t_{2n})])
\end{eqnarray*}
defines a class in the Hopf-invariant cyclic cohomology, $\text{ch}^{2n}_{\Theta} \in HC^{2n}_{\Theta}(A)$.
\end{quote}

\noindent Proof. The pairing of item 1 is independent on the choice of the K-theory class and cyclic cohomology class which is proved as in the nonequivariant case. If the  summability of $D_q$ on $\hil$ is $2n$, then it is $2n$ or less in $\hil^{\Theta}$. Then the commutators $[F_q, \pi(t_{i})]$ are in the Schatten ideal $\mathcal{L}^{2n+1}(\hil^{\Theta})$ and so by Holder inequality, $\text{ch}^{2n}_{\Theta}$ is well defined. It is well known that the functional of item 2 is in the kernel of $b$ and so defines a cocycle. \5 $\square$\3

Let $p$ be a representative of a Hopf-equivariant projective module in $\text{end}(V) \otimes \mathbb{C}[G_q]$. The primitive coproduct in $\mathbb{C}[K]$ makes  $V \otimes \hil$ a $\mathbb{C}[K]$-module and we can take the invariant subspace. The projections
\begin{eqnarray*}
p_{\pm} = \frac{1 \pm \omega}{2} \pi(p)  
\end{eqnarray*}
commute with the action of $\mathbb{C}[K]$ and so they restrict to projection operators on the invariant subspace of $V \otimes \hil$. Denote by $\hil^{\Theta \pm}_p$ the $\mathbb{C}[K]$-invariant ranges. The operator $p D_q p$ is $\mathbb{C}[K]$-invariant and off diagonal with respect to the chiral grading, and so we have the components 
\begin{eqnarray*}
D_{pq}^+ = p_- D_q p_+ : \hil_p^{\Theta +} \rightarrow \hil_p^{\Theta -} \5 \text{and} \5 D_{pq}^- =  p_+ D_q p_-: \hil_p^{\Theta -} \rightarrow \hil_p^{\Theta +}. 
\end{eqnarray*}
As usual we can define the approximate sign of $D_{pq}$ which is given by 
\begin{eqnarray*}
F_{pq}^{\pm} = p_{\mp} \frac{D_{q}}{\sqrt{1 + D_{q}^2}} p_{\pm}. 
\end{eqnarray*}
The index of the component $F_{pq}^+$ depends only on the Hopf-equivariant K-theory class of $p$. One can compute the index by pairing the character $\text{ch}_{\Theta}^{2n}$, with $n$ sufficiently large, with the K-theory element $[p]$. The next theorem is easily verified with the standard analysis of Fredholm index problems in noncommutative geometry. 

\3 \noindent \textbf{Theorem 3.} If $2n$ is the dimension of $G$, then the index of $F_{pq}^+$ on the invariant subspace $\hil^{\Theta +}_p$ is given by the pairing
\begin{eqnarray*}
\text{ind}(F_{pq}^+) = ([\text{ch}_{\Theta}^{2n}], [p])
\end{eqnarray*}
for any Hopf-equivariant projective module $p$.

\3 \noindent Proof. If $P: \hil^1 \rightarrow \hil^2$ is a Fredholm operator and $Q: \hil^2 \rightarrow \hil^1$ its parametrix ($\hil_1$ and $\hil_2$ are complex Hilbert spaces) and if $1 - PQ$ and $1-QP$ are both in the Schatten ideal $\mathcal{L}^{\frac{2n+1}{2}}$, then the  index of $P$ can be computed by \cite{Cal67}, \cite{Con85}
\begin{eqnarray*}
\text{ind}(P) = \Tr ( (1-PQ)^k - (1-QP)^k )
\end{eqnarray*}
for any integer $k \geq (2n+1)/2$. Following \cite{Con85}, we write $P = F^+_{pq} : \hil^{\Theta +}_p \rightarrow \hil^{\Theta -}_p$ and $Q = F_{pq}^-: \hil^{\Theta -}_p \rightarrow \hil^{\Theta +}_p$. Then the operators $1 - PQ$ and $1 - QP$ are the $++$ and $--$ components of 
\begin{eqnarray*}
p - pF_q pF_q p = - p [F_q,p]^2 p \in  \mathcal{L}^{\frac{2n + 1}{2}}.
\end{eqnarray*}
A standard algebraic manipulation leads to
\begin{eqnarray*}
\frac{(-1)^n}{2} \Tr (\omega F_q[F_q,p]^{2n + 1} ) = \Tr( \omega(p - pF_q pF_q p)^{n + 1})
\end{eqnarray*}
and the right hand side computes the index because $\omega = \text{diag}(1,-1)$. \5 $\square$ \3

\noindent \textbf{4.5.} Let $\Theta \rightrightarrows G_q$ be a quantum orbifold associated with a finite subgroup $K$ in $T$. If the action $\Theta$ of $K$ on $\mathbb{C}[G_q]$ is faithful, then $\Theta \rightrightarrows G_q$ is defined to be an effective quantum orbifold.  Consider the crossed product algebra $\mathbb{C}[K \ltimes G_q]$ of 2.1. This algebra can be represented on the space $\hil$ by setting
\begin{eqnarray}\label{44}
\varpi(\sigma \cdot t) \psi = \pi(t)(\sigma \triangleright \psi)  
\end{eqnarray}
for all $\sigma \in K$, $t \in \mathbb{C}[K \ltimes G_q]$ and $\psi \in \hil$, and extended linearly.

\3 \noindent \textbf{Theorem 4.} Let $\Theta \rightrightarrows G_q$ be an effective spin quantum orbifold associated with a finite group $K$. Then the collection $(\mathbb{C}[K \ltimes G_q], D_q, \hil)$ defines a dim($G$)-summable spectral triple which is even if dim($G$) is even. Up to a unitary equivalence of spectral triples, this model is independent on the choices $\mathcal{F}$ and $\phi$. 

\3 \noindent Proof. Denote by $\sigma_i$, $i = \{1, \ldots, p\}$ the set of group elements. To see that $\varpi$ is a faithful representation it is sufficient to check that the operators $\varpi(\sigma_i \cdot t^{\lambda}_{\mu \nu})$ are linearly independent for all possible indices $i, \mu, \nu$ and $\lambda$. Suppose that $\varpi(\sigma_i \cdot t^{\lambda}_{\mu \nu})$ and $\varpi(\sigma_{i'} \cdot t^{\lambda'}_{\mu' \nu'})$ are linearly dependent and let these operators act on the subspace $1 \otimes \Sigma$.
 Now we see that $t^{\lambda}_{\mu \nu} = t^{\lambda'}_{\mu' \nu'}$ must hold and $\hat{\varrho}(\sigma_i) = \hat{\varrho}(\sigma_j)$ must hold. Then we fix $v \in \Sigma$ and let the operators act on $\mathbb{C}[G_q] \otimes v$. It then follows from the faithfulness of the action $K$ on $\mathbb{C}[G_q]$ that $\sigma_i = \sigma_j$ must hold. So, the representation $\varpi$ is faithful. 

Since $D_q$ and $\hil$ are given exactly as in the quantum group model which is an isospectral model, the spectral properties of the Dirac operator are as in the case of the Lie group $G$. In particular, the summability is as claimed. The $K$-invariance of $D_q$ gives us the following commutator formula
\begin{eqnarray*}
[D_q, \varpi(\sigma \cdot t)] \psi &=& D_q \pi(t)(\sigma \triangleright \psi) - \pi(t) (\sigma \triangleright (D_q \psi)) \\
&=& [D_q, \pi(t)] (\sigma \triangleright \psi).
\end{eqnarray*}
if $\psi$ is in the domain of $D_q$. It follows that the commutators $[D_q, \varpi(\sigma \cdot t)]$ extend to bounded operators on $\hil$. In the even dimensional cases, the quantum group model has the chirality operator $\omega$ which anticommutes with $D_q$ and $\omega$ commutes with the representation $\varpi$ because the action of $K$ and the representation $\pi$ commute with the chirality operator.

The choice of the twist does not have an effect on the Dirac operator. In the proof of Theorem 2 we observed that if $\phi$ and $\phi'$ is a pair of $*$-homomorphism $\mathcal{U}(G_q) \rightarrow \mathcal{U}(G)$, then the resulting Dirac operator are unitary conjugates, and the action of $K$ commutes with the unitary operator. So, the unitary operator commutes with the representation $\varpi$ and so induces a unitary equivalence of spectral triples. \5 $\square$\3 

It would be  interesting to understand the structure of the crossed product quantum orbifold $C^*$-algebras, in particular their K-theory. On a classical topological level, the equivariant K-theory classifies equivariant vector bundles and the K-theory groups are  very accurate tools to capture information of the isotropy of the action. If the action is free, then the equivariant K-theory is isomorphic to the ordinary K-theory of the quotient. The difference between the equivariant K-theory and the K-theory of the quotient measures the singularity of the action. In the quantum orbifold model, the problem if finding singularities can be approached by studying the difference of the K-theory groups of the crossed product $C^*$-algebra $C(K \ltimes G_q)$ and the fixed point $C^*$-algebra $C(G_q)^{\Theta}$. For example, consider the quantum teardrop corresponding to the parameters $(k,l) = (1,l)$ introduced in 2.2. Classically such teardrop orbifold would have a $\mathbb{Z}_l$ isotropy but being homeomorphic to $S^2$, the quotient space would be insensitive to the effects of the isotropy group. On the other hand, the K-theory groups for the fixed point quantum teardrop $C^*$-algebras are given by $K_0(C(G_q)^{\Theta}) \simeq \mathbb{Z}^{l+1}$ and $K_1(C(G_q)^{\Theta}) \simeq 0$. So, it would seem that this model is sensitive to the effects of the isotropy already on the quotient level. However, the K-theory of the fixed point algebra is insensitive to the change the parameter $k$, and for this reason, one should study the algebra $C(K \ltimes G_q)$ if $k > 1$. This discussion is also in agreement with the results of \cite{Bre14}.

\end{document}